\RequirePackage{fix-cm}
\documentclass[smallcondensed]{svjour3}
\usepackage{amssymb,amsmath}
\usepackage{url}
\usepackage[hidelinks]{hyperref}
\usepackage{graphicx}

\newcommand{\Z}{\mathbb{Z}}
\newcommand{\F}{\mathbb{F}}
\begin{document}

\title{A new distance-regular graph of diameter $3$\\ on $1024$ vertices%
\thanks{This research is supported by National Natural Science Foundation of China (61672036),
Technology Foundation for Selected Overseas Chinese Scholar, Ministry of Personnel of China (05015133), Excellent Youth Foundation of Natural Science Foundation of Anhui Province (1808085J20) and
Key projects of support program for outstanding young talents in Colleges and Universities (gxyqZD2016008).}
}


\author{Minjia~Shi 
\and 
\href{https://orcid.org/0000-0002-8516-755X}{Denis~S.~Krotov}
\and
\href{https://orcid.org/0000-0002-4078-8301}{Patrick~Sol\'e}
}

\authorrunning{M.~Shi, D.~Krotov, P.~Sol\'e} 

\institute{M.~Shi \at
Key Laboratory of Intelligent Computing Signal Processing, 
Ministry of Education, Anhui University, 
No.3 Feixi Road, Hefei, Anhui, 230039, China; \\
School of Mathematical Sciences, Anhui University, Hefei, Anhui, 230601,
China \\
              \email{smjwcl.good@163.com}           
             \and
            D.~S.~Krotov \at
            Sobolev Institute of Mathematics, 
            pr. Akademika Koptyuga 4, 
            Novosibirsk 630090, Russia \\
            \email{krotov@math.nsc.ru}
                          \and
             P.~Sol\'e \at
             4CNRS/LAGA, University of Paris 8, 
             2 rue de la Libert\'e, 
             93 526 Saint-Denis, France \\
             \email{sole@enst.fr}
}

\date{Received: 2018-06 / Accepted: not yet }

\maketitle

\begin{abstract}
The dodecacode is a nonlinear additive quaternary code of length $12$. By puncturing it at any of the twelve coordinates, we obtain a  uniformly packed code of distance $5$.
In particular, this latter code is completely regular but not completely transitive.
Its coset graph is distance-regular of diameter three on $2^{10}$ vertices,
with new intersection array $\{33,30,15;1,2,15\}$.
The automorphism groups of the code, and of the graph, are determined. Connecting the vertices at distance two
gives a strongly regular graph of (previously known) parameters $(2^{10}, 495,238, 240)$. Another strongly regular graph with the same parameters is constructed on the codewords of the dual code.
A non trivial completely regular binary code of length $33$ is constructed.
\end{abstract}
\section{Introduction}\label{s:intro}
Distance-regular (DR) graphs form the most extensively studied class of structured graphs due to their many connections with codes, designs, groups and orthogonal polynomials \cite{BanIt-1,Brouwer}. Since the times of Delsarte \cite{Delsarte:1973}, a powerful way to create DR graphs, especially in low diameters, has been to use the coset graph of completely regular codes. A code is {\em completely regular} if the weight distribution of each coset solely depends on the weight of its coset leader. In \cite{Brouwer} can be found many such examples from Golay codes, Kasami codes, and others. A recent survey is \cite{BRZ:CR}.
A special class of completely regular codes is that of uniformly packed codes. A code with packing radius $e$ is {\em uniformly packed} if the spheres of radius $e+1$ about the codewords overlap in a very uniform way: there are two constants $\lambda$ and $\mu$ such that the number of codewords at distance $e+1$ from some $x$ in ambient space is either $\lambda$ if $x$ is at distance $e$ from the code or $\mu$ if $x$ is at distance $\ge e+1$ from the code.
In 1974 the following hypothetical parameters $[n,k,\ge 2e+1]$ for quaternary uniformly packed codes with packing constants $(\lambda,\mu)$ were found \cite{BZZ:1974:UPC} by computer search as
$$n=\frac{2^{2m+1}+1}{3},\quad n-k=2m+1, \quad e=2,\quad\mu=\lambda+1=\frac{2^{2m}-1}{3}$$
for $m\ge 2$, but no corresponding code was found. The case $m=2$ of the above parameters leads to a putative coset graph listed in \cite[p.428]{Brouwer} as a distance-regular graph on $1024=4^5$ vertices of diameter $3$. It does not appear to be solved as of 2016 in the recent survey \cite{vDKT:DR}.
In 1998, an additive quaternary code of parameters $(12,4^6,6)$ was introduced in \cite{CRSS:quantum} for the purpose of quantum-error correction, and called the {\em dodecacode} . This code was used in \cite{KimPless:GF4} to construct designs by shortening.
In this note, we construct an additive non-linear code with the Bassalygo et al. parameters for $m=2$ by puncturing the dodecacode. From the coset graph of that code,
we obtain a first known  distance-regular graph with intersection array $ \{33,30,15;1,2,15\}$, already studied in \cite{Makh:33-30-15-1-2-15}, as a hypothetical object,  from the standpoint of symmetry. This is the main result of this note. As a byproduct, connecting the vertices at distance two
gives a strongly regular graph of parameters $(2^{10}, 495,238, 240)$. Surprisingly, connecting the codewords in the dual code at distance $8$ gives also a strongly regular graph with the same parameters.
Other results coming from quaternary additive codes are as follows.
By triple puncturing of the dodecacode, we obtain an additive code of length $9$, and taking its coset graph,  we construct a strongly regular graph of parameters $(64,27,10,12)$. We give a universal correspondence between quaternary and binary codes that preserves coset graphs, up to isomorphism.
When applied to the two preceding codes, it constructs new completely regular binary codes in length $27$ and $33$.
  Thus, the observation that the punctured dodecacode is uniformly packed solves a problem that had been open since 1974 for quaternary codes and since 1989 for distance-regular graphs.

 The sections are arranged as follows. The next section collects the necessary notation and definitions. Section~\ref{s:codes} studies codes either quaternary or binary. Section~\ref{s:graphs} constructs the new distance regular graph announced in the title. Section~\ref{s:concl} concludes the article and points out some challenging open problems.

 As finding a new distance-regular graph is a most exciting result of the current work, we give here one of its shortest descriptions: it is the Cayley graph on $\Z_2^{10}$ with the connecting set
\{1, 2, 4, 8, 16, 32, 54, 64, 128, 149, 151, 170, 186, 216, 217, 256, 293, 310, 329, 338, 466, 512, 597, 605, 658, 681, 745, 841, 951, 952, 956, 966, 998\},
where integers are treated as binary tuples in the standard way.
\section{Definitions and notation}\label{s:def}
For any undefined term pertaining to codes we refer to \cite{HuffmanPless,MWS}, and to \cite{BanIt-1,Brouwer} for any undefined term related to distance-regular graphs.
\subsection{Linear codes}\label{ss:linear}
A  linear code of length $n$, dimension $k$, minimum distance $d$ is called an $[n,k,d]$ code. A two-weight code is an $[n,k]$ linear binary code having two nonzero weights $w_1$ and $w_2$. Its parameters are denoted by the formula
$[n,k;w_1,w_2]$. The duality is understood with respect to the standard inner product. The {\em external distance} of a linear code is the number of nonzero weights in its dual.
A {\em coset} of a linear code $C$ is any translate of $C$ by a constant vector. A {\em coset leader} is any coset element that minimizes the weight. The {\em weight of a coset} is the weight of any of its leaders.
The {\em coset  graph} $\Gamma_C$ of a code $C$ is defined on the cosets of $C$, two cosets being connected if they differ by a coset of weight one.
\subsection{Additive codes}\label{ss:additive}
We consider codes of length $n$ for the Hamming distance over $\F_4$, and denote the Hamming weight of $x$ by $\mathrm{wt}(x)$.
An additive code of length $n$ over $\F_4$ is an additive subgroup of $\F_4^n$. It is a free $\F_2$ module with $4^k$ elements for some $k\le n$ (here $2k$ is an integer, but $k$ may be half-integral). There are two ways to specify such a
code $C$. Firstly, using a {\em generator matrix} $G$, the code $C$ can be cast as the $\F_2$-span of its rows. Secondly, a {\em parity check matrix} $H$ can be defined by using the {\em trace inner product} given
in dimension $1$ by $x*y=\mathrm{Tr}(xy^2)$, where $\mathrm{Tr}(z)=z+z^2$, and extended coordinatewise to vectors of length $n$. The {\em dual} of a code $C$ denoted by $C^\bot$ is understood with respect to that inner product.
Define the $\star$-product of a $m\times n$  matrix $M$ with rows $M_i$ by a vector of length $n$ as
$$ M\star x= (M_1 *x,\cdots, M_m*x)^T, $$ the $^T$ denoting transposition.
An additive code with $4^k$ codewords can then be specified by an $2n-2k$ by $n$ matrix $H$, the {\em parity check matrix} of $C$, as
$$ C=\{x \in  \F_4^n \mid H\star x=0\}. $$
Note that $H$ is a generator matrix for $C^\perp$.
 We call {\em syndrome} $s(x)$ of $x\in \F_4^n$ the column vector of size $2n-2k$ defined by $s(x)=H\star x$.
A code $C^-$ is obtained by {\em puncturing} from a code $C$ at the coordinate $i$ if it is the projection of $C$ on the remaining $n-1$ other coordinates.
A {\em coset} of an additive code $C$ is any translate of $C$ by a constant vector. The definitions of coset leader and coset weight are the same as in the case of linear codes.
The {\em coset  graph} $\Gamma_C$ of a code $C$ is then the graph defined on the $4^{n-k}$ syndromes, two of them being connected if they differ by a syndrome $s(x)$ with $\mathrm{wt}(x)=1$.
By the obvious one-to-one correspondence between syndromes and cosets \cite{HuffmanPless}, this graph is also the graph on the cosets, two cosets being connected if they differ by a coset of weight one.
We give without proof the following extension of \cite[11.1.11]{Brouwer} from linear to additive codes.
{\theorem \label{ext} If $C$ is an additive quaternary code of minimum distance at least three, with dual weight distribution $[\langle i,A_i\rangle]$, then the spectrum of $\Gamma_C$ is $\{(3n-4i)^{A_i}\}$. Thus $A_i$ is the frequency
of weight $i$ in $C^\bot$ and the multiplicity of the eigenvalue $3n-4i$.}

The {\em monomial automorphism group} of a quaternary code $C$ is the set of all monomial transformations that leave the code wholly
invariant, where a {\em monomial transformation} \cite{HuffmanPless}
is the composition of
a coordinate permutation and the multiplication of the value in every position by a nonzero scalar. The \emph{complete automorphism group} $\mathrm{Aut}(C)$ of a code $C$ is its stabilizer in the group of isometries of the Hamming space.
The complete automorphism group includes the monomial automorphism group and the group of translations of the code, but can be larger than their product.
An additive quaternary code is called {\em completely transitive} in the sense of \cite{Sole:CR&CTC} if the induced action of $\mathrm{Aut}(C)$ on the cosets of given weight of $C$ is transitive. Such codes are completely regular, with a distance-transitive coset graph. For an alternative definition see \cite{GiuPra:1999}. For a partial classification see \cite{BRZ:2001:e3}.
\subsection{Graphs}
All graphs in this note are finite, undirected, connected, without loops or multiple edges. The neighborhood $\Gamma(x)$ is the set of vertices connected to $x$.
The {\em degree} of a vertex $x$ is the size of $\Gamma(x)$.
A graph is {\em  regular} if every vertex has the same degree. The $i$-neighborhood $\Gamma_i(x)$ is the set of vertices at geodetic distance $i$ to $x$.
A graph is {\em  distance regular} (DR) if for every pair or vertices $u$ and $v$ at distance $i$ apart the quantities
\begin{eqnarray*}
 a_i&=&| \Gamma_{i+1}(u)\cap \Gamma(v)|\\
 c_i&=&| \Gamma_{i-1}(u)\cap \Gamma(v)|
\end{eqnarray*}
solely depend on $i$ and not on the special choice of the pair $(u,v)$.
In that situation the graphs $\Gamma_i$ are regular of degree $v_i$ and we will refer to the $v_i$s as the {\em successive degrees} of $\Gamma$.
The {\em  automorphism group} of a graph is the set of permutations of the vertices that preserve adjacency.
A graph is {\em distance-transitive} if its automorphism group is transitive on its vertices and on each of the sets $\Gamma_i(x)$ for any $i$ and any fixed $x$.
A DR graph of diameter $2$ is called a {\em Strongly Regular Graph} (SRG). Its parameters $(\nu,\kappa,\lambda,\mu)$ denote the number of vertices, the degree, the number of common neighbors of a pair of connected vertices, the number
 of common neighbors of a pair of disconnected vertices.
 The {\em spectrum} of a graph is the set of distinct eigenvalues of its adjacency matrix. It is denoted by $\{\lambda_1^{m_1},\lambda_2^{m_2},\dots\}$, where $m_i$ stands for the multiplicity of the eigenvalue $\lambda_i$.
 The {\em Hamming graph} $H(n,q)$ is the DR graph on $\F_q^n$, two vectors being connected if they are at Hamming distance one.

\section{Codes}\label{s:codes}
\subsection{Quaternary codes}\label{ss:quaternary}
The dodecacode $D$ is a code of length $12$ over $\F_4=\F_2(w)$ of parameters $(12,4^6,6)$ that is additive but not linear \cite{CRSS:quantum}. It can be defined as a cyclic code with one generator. In the notation of \cite{CRSS:quantum} we have
 $D=(w10100100101)$.
Thus, puncturing at any of the twelve coordinates give an equivalent code $D^-$ of parameters  $(11,4^6,5)$, of generator matrix given in Table 1.
\begin{center}
\small {Table 1: generator matrix of $D^-$} \vspace{-0.8cm}
\end{center}
\begin{center}
\small
\begin{eqnarray*}
                                               \left(
      \begin{array}{ccccccccccc}
        0 &  0 &  0 &  0 &  0& w^2& w^2&   0 &  w &  1&   w\\
    0&   0 &  0&   0&   0&   w &  0&   w&   w&   w&   1\\
     1&   0&   0&   0 &  0 &  1&   0 &  1& w^2& w^2&   1\\
    w &  0 &  0 &  0 &  0 &  0 &  w &  1&   w&w&   w\\
     0&   1&   0 &  0&   0&   0 &  1  & 1&  1& w^2& w^2\\
    0 &  w &  0 &  0&   0& w^2&   1& w^2&   w&   w&   0\\
   0 &  0   &1 &  0 &  0 &  w &  w &  1 &  1&  0&   1\\
     0 &  0 &  w &  0 &  0&   w &  1&   w& w^2& w^2&   0\\
     0 &  0&   0&   1&   0 &  w&   w&   w&   0&  1 &  w\\
     0 &  0 &  0 &  w &  0& w^2 &  1 &  1& w^2&   0&   w\\
    0&   0&   0&   0&   1 &w^2& w^2& w^2&   1&   0 &w^2\\
    0 &  0&   0 &  0 &  w&   w &  1&   0 &  1&   w&   w\\
      \end{array}
    \right)
    \end{eqnarray*}

\end{center}
 Its primal and dual weight distributions and automorphism groups are easily computed in Magma \cite{Magma} and Sage~\cite{sage}.
{\theorem \label{distri}\label{prim} The  weight distribution and the dual weight distribution of $D^-$ are
\begin{eqnarray*}
&[ \langle 0, 1\rangle , \langle 5, 198\rangle , \langle 6, 198\rangle , \langle 7, 990\rangle , \langle 8, 495\rangle , \langle 9, 1650\rangle , \langle 10, 330\rangle , \langle 11,
234\rangle  ] \\ &\mbox{and} \quad[ \langle 0, 1\rangle , \langle 6, 198\rangle , \langle 8, 495\rangle , \langle 10, 330\rangle  ], 
\end{eqnarray*}
respectively, where $\langle i,A_i\rangle $ means that there are $A_i$ codewords of weight $i$.}

{\theorem The monomial automorphism group of $D^-$ is generated by the following two monomial matrices:
$$
\left[\begin{array}{c@{\ \,}c@{\ \,}c@{\ \,}c@{\ \,}c@{\ \,}c@{\ \,}c@{\ \,}c@{\ \,}c@{\ \,}c@{\ \,}c}
0&0&0&0&0&\makebox[0mm]{$w$}&0&0&0&0&0\\
1&0&0&0&0&0&0&0&0&0&0\\
0&0&\makebox[0mm]{$w^2$}&0&0&0&0&0&0&0&0\\
0&0&0&0&0&0&0&\makebox[0mm]{$w^2$}&0&0&0\\
0&0&0&0&0&0&0&0&0&1&0\\
0&0&0&0&\makebox[0mm]{$w$}&0&0&0&0&0&0\\
0&0&0&0&0&0&0&0&0&0&\makebox[0mm]{$w$}\\
0&0&0&1&0&0&0&0&0&0&0\\
0&\makebox[0mm]{$w^2$}&0&0&0&0&0&0&0&0&0\\
0&0&0&0&0&0&0&0&\makebox[0mm]{$w^2$}&0&0\\
0&0&0&0&0&0&1&0&0&0&0\\
\end{array}\right],\qquad
\left[\begin{array}{c@{\ \,}c@{\ \,}c@{\ \,}c@{\ \,}c@{\ \,}c@{\ \,}c@{\ \,}c@{\ \,}c@{\ \,}c@{\ \,}c}
1&0&0&0&0&0&0&0&0&0&0\\
0&0&0&0&0&0&0&0&0&0&\makebox[0mm]{$w^2$}\\
0&0&0&0&0&\makebox[0mm]{$w^2$}&0&0&0&0&0\\
0&0&0&0&0&0&0&\makebox[0mm]{$w^2$}&0&0&0\\
0&0&0&0&0&0&0&0&\makebox[0mm]{$w$}&0&0\\
0&0&\makebox[0mm]{$w$}&0&0&0&0&0&0&0&0\\
0&0&0&0&0&0&0&0&0&1&0\\
0&0&0&\makebox[0mm]{$w$}&0&0&0&0&0&0&0\\
0&0&0&0&\makebox[0mm]{$w^2$}&0&0&0&0&0&0\\
0&0&0&0&0&0&1&0&0&0&0\\
0&\makebox[0mm]{$w$}&0&0&0&0&0&0&0&0&0\\
\end{array}\right].
$$
It is isomorphic to a subgroup of $S_{33}$ of order $54=2 \times 3^3$ with generators 
    $(1$ $4$ $27$ $29$ $14$ $18)(2$ $5$ $25$ $30$ $15$ $16)(3$ $6$ $26$ $28$ $13$ $17)(7$ $9$ $8)(10$ $22$ $12$ $24$ $11$ $23)(19$ $31$ $20$ $32$ $21$ $33)$ and
   $(4$ $32)(5$ $33)(6$ $31)(7$ $17)(8$ $18)(9$ $16)(10$ $23)(11$ $24)(12$ $22) (13,27)(14$ $25)(15$ $26)(19$ $28)(20$ $29)(21$ $30)$ in disjoint cycles product notation (this subgroup reflects the action of the automorphism group on the words of weight $1$ in the order $10...0$, $w0...0$, $w^20...0$, $010...0$, \ldots, $0...0w^2$).
The complete automorphism group of $D^-$ is the product of the monomial automorphism group with the group of translations and has the order $54\cdot 2^{12}$.}

\begin{remark}
Since $54$ is not divisible by $33$ the code $D^-$ cannot be completely transitive \cite{Sole:CR&CTC}.
\end{remark}

        To compute the parameters $\lambda$ and $\mu$, we require the following lemma, of independent interest.
        {\lemma If $C$ is a uniformly packed code of distance $2e+1$, length $n$, and weight distribution $A_i$ over $\F_q$, then we have
        \newcommand\bnm[2]{\bigg(\begin{array}{@{}c@{}}#1\\#2\end{array}\bigg)}
        \begin{eqnarray*}
        \lambda (q-1)^e{\bnm ne} &=&A_{2e+1}\bnm{2e+1}{e}\\
        A_{2e+1}\bnm{2e+1}{e}(e+1)(q-2)+A_{2e+2}\bnm{2e+2}{e+1}&=&(\lambda-\mu)A_{2e+1}\bnm{2e+1}{e}\\
        &&{}+(\mu-1)(q-1)^{e+1}\bnm{n}{e+1}.
        \end{eqnarray*}}\vspace{-0.5cm}
        \begin{proof}
        Double counting the number of pairs $(x,y)$ such that 
        $\mathrm{wt}(x)=e$, 
        $\mathrm{wt}(y)=2e+1$, 
        $\mathrm{wt}(y-x)=e+1$,
        $y\in C$ relates
$\lambda$ and $A_{2e+1}$. This yields the first relation.


Double counting the number of pairs $(x,y)$ such that $\mathrm{wt}(x)=e+1$, 
$\mathrm{wt}(y)\in \{2e+1,2e+2\}$, 
$\mathrm{wt}(y-x)=e+1$,
$y\in C$
(note that among the $(q-1)^{e+1} \binom{n}{e+1} $ words of weight $e+1$, exactly $A_{2e+1} \binom{2e+1}{e}$ are at distance $e$ from the code,
and the other are at distance $e+1$ from the code)
relates
$\lambda$, $\mu$, $A_{2e+1}$, and $A_{2e+2}$. This gives us the second relation.
        \qed\end{proof}
        
        We are now in a position to state and prove the main result of this subsection.
        {\theorem \label{lambmu} The code $D^-$ is uniformly packed with $(\lambda,\mu)=(4,5)$.}
        
        \begin{proof}
We know the code is uniformly packed since its minimum distance equals twice its external distance minus one \cite[Cor. 11.1.2]{Brouwer}. To compute $\lambda$ and $\mu$ we specialize the above lemma to $q=4, \,n=11,\, e=2$,
and $A_5=A_6=198$ that is known from Theorem~\ref{prim}.
The first relation yields $\lambda=4$. Reporting into the second gives $\mu=5$.
        \qed\end{proof}

        Denote by $D^{3-}$ any $(9, 4^6)$ code obtained by puncturing $D^-$ on any pair of coordinates.  The dual of $D^{3-}$ is an additive two-weight code.
        The following result comes from an easy Magma computation \cite{Magma}.
        {\theorem \label{distri2} The dual weight distribution of $D^{3-}$ is $[ \langle 0, 1\rangle , \langle 6, 36\rangle , \langle 8, 27\rangle  ]$.}

         \subsection{Binary codes}\label{ss:binary}
        We aim to construct a binary code with a coset graph isomorphic to that of $D^-$. To that end, we define a universal correspondence between quaternary and binary codes that preserves coset graphs, up to isomorphism.

        Define the concatenation map of a quaternary code of length $n$ with the zero sum code 
        $R_3^\bot=\{000,011,110,101\}$ 
        as follows.
        For binary scalars $a$ and $b$, let
        $\phi(a+bw)=(b,b+a,a)$. Equivalently,
        $\phi(c)=(\mathrm{Tr}(c),\mathrm{Tr}(wc),\mathrm{Tr}(w^2c))$. The action of $\phi$ is extended to the quaternary vectors coordinatewise.
        
        
        {\theorem 
        \label{th:qua-to-bin}
        Assume that $Q$ is a quaternary additive code of length $n$ and $B=\phi(Q^\bot)^\bot$. The following assertions hold:
        \begin{itemize}
         \item[\rm(i)] if $Q$ is an $(n,4^k)$ code, then $B$ is a linear $[3n,n+2k]$ binary code;
        \item[\rm(ii)]         if $w$ is a weight of $Q^\perp$ with frequency $A_w$, then $2w$ is a weight of $B^\bot$ with frequency $A_w$, and all the weights of $B^\bot$ arise in this way;
         \item[\rm(iii)] the coset graphs $\Gamma_Q$ and $\Gamma_B$ are isomorphic.
         \end{itemize}

        }
\begin{proof} 
(i) Trivially, the length of $B^\bot = \phi(Q^\bot)$ is $3n$ and the size is the same as the size of $Q^\bot$, i.e., $4^{n-k}$. 
Hence, the size of $B$ is $2^{3n}/4^{n-k}=2^{n+2k}$. 

(ii) From the concatenation, we readily see
that $\mathrm{wt}(\phi(c))=2\mathrm{wt}(c)$.

(iii) We will show that the graphs $\Gamma_Q$ and $\Gamma_B$
built on the syndromes of the check matrices $P$ and $\phi(P)$ coincide.
To see this, we first note that the both graphs are built on the binary columns of height $2n-2k$.
Next, we consider the connecting syndromes that correspond to the weight-$1$ vectors.
Denote by $P_1$ the first column of the check matrix. The syndromes corresponding to the 
weight-$1$ quaternary vectors $(1,0,...,0)$, $(w,0,...,0)$, and $(w^2,0,...,0)$
are $\mathrm{Tr}(P_1)$,  $\mathrm{Tr}(w^2 P_1)$, and $\mathrm{Tr}(w P_1)$, where the trace map acts on the column component-wise. On the other hand, by the definition of $\phi$,
the first three columns of $\phi(P)$ are $\mathrm{Tr}(P_1)$,  $\mathrm{Tr}(w P_1)$, and $\mathrm{Tr}(w^2 P_1)$. Therefore, the syndromes corresponding to the binary $1$-weight vectors
$(1,0,...,0)$, $(0,1,0,...,0)$, $(0,0,1,0,...,0)$ are $\mathrm{Tr}(P_1)$,  $\mathrm{Tr}(w P_1)$,  $\mathrm{Tr}(w^2 P_1)$ again. Considering in a similar way every other column of $P$, 
we find that the set of syndromes corresponding  to the weight-$1$ vectors is the same for $P$ and $\phi(P)$. Hence, $\Gamma_Q = \Gamma_B$.
\qed\end{proof}

\begin{remark}
 It can be shown that $Q$ and $B$ are related by the following correspondence (denoted by $\psi$) that associates $2^n$ codewords of $B$ to any codeword of $Q$. For simplicity, we write it for $n=1$.
\begin{eqnarray*}
0   \stackrel{\psi}\longrightarrow 000,111; \qquad
1   \stackrel{\psi}\longrightarrow 100,011;\qquad
w   \stackrel{\psi}\longrightarrow 001,110;\qquad 
w^2 \stackrel{\psi}\longrightarrow 010,101.
\end{eqnarray*}
Thus the images of any of the four symbols form an antipodal pair of vertices in the $3$-cube.
We have the following commutative diagram where the down arrow means ``dual''.
\[\begin{array}{ccc} Q &
{\longrightarrow} &
\psi(Q)\\
\big\downarrow & &
\big\downarrow\vcenter{%
\rlap{}}\\
Q^{\bot} & \longrightarrow & \phi(Q^{\bot})=\psi(Q)^\bot
\end{array}\]
Thus $\psi$ is the pullback of $\phi$ in this diagram.
\end{remark}

        {\corollary Let $B^-=\phi({D^-}^\bot)^\bot$. The binary code $B^-$ is a completely regular code of parameters $[33,23,3]$, with dual weight distribution 
        $$
        [ \langle 0, 1\rangle , \langle 12, 198\rangle , \langle 16, 495\rangle , \langle 20, 330\rangle  ], $$ 
        and weight distribution
        \begin{multline}\nonumber
        [ \langle 0, 1\rangle , \langle 3, 11\rangle , \langle 5, 198\rangle , \langle 6, 1243\rangle , \langle 7, 4158\rangle , \langle 8, 13563\rangle , 
        \langle 9, 38445\rangle , \\
        \langle 10,88638\rangle , \langle 11, 185397\rangle , \langle 12, 352902\rangle , \langle 13, 568788\rangle , 
        \langle 14, 786885\rangle , \qquad\\  \langle 15, 998052\rangle ,
\langle 16, 1156023\rangle ,
\langle 17, 1156023\rangle , \langle 18, 998052\rangle , 
\langle 19, 786885\rangle ,\\\qquad
\langle 20, 568788\rangle ,
\langle 21,352902\rangle , \langle 22, 185397\rangle , \langle 23, 88638\rangle , \langle 24, 38445\rangle , 
\\
\langle 25, 13563\rangle , \langle 26, 4158\rangle , \langle 27,1243\rangle , \langle 28, 198\rangle , \langle 30, 11\rangle , \langle 33, 1\rangle  ]. 
\end{multline}
        } \vspace{-0.8cm}
        \begin{proof} The dual weight distribution follows by the preceding theorem combined with Theorem \ref{distri}. The weight distribution is then computed by MacWilliams transform.
        \qed\end{proof}

\begin{remark}
The code $B^-$ is a non trivial example of a completely regular code since its minimum distance is $3$, and its external distance is $3$. Thus it is neither perfect, nor uniformly packed.
\end{remark}

    {\corollary The code $\phi({D^{3-}}^\bot)$ is a binary two-weight $[27,6;12,16]$ code.}
    \begin{proof} The  weight distribution follows by Theorem~\ref{th:qua-to-bin}, where $n=9$, $k=6$, $Q=D^{3-}$, and $B=\phi({D^{3-}}^\bot)^\bot$, combined with Theorem \ref{distri2}.
        \qed\end{proof}
        \begin{remark}
         Binary completely regular codes with the last parameters are known \cite{CRCtable2}; all such codes are related to bent functions, see Theorem~12.12 and the following paragraph in \cite{CalKan:86:2weight}. Puncturing different
         coordinates, we obtain $3$ nonequivalent codes with these
         parameters.
        \end{remark}
    \section{Graphs}\label{s:graphs}
    In this section, we study the coset graphs of ${D^-}$, and of $D^{3-}$.
    {\theorem The graph $\Gamma_{D^-}$ is distance-regular of diameter $3$, of spectrum $$\{27^1,9^{198},1^{495},(-7)^{330}\}. $$ Its successive degrees are $(1,33,495,495)$.}
    \begin{proof}
    The spectrum of $\Gamma_{D^-}$ is easily computed from the weight distribution of Theorem \ref{distri} upon applying Theorem \ref{ext}. The distance-regularity follows from the fact that $D^-$ is completely regular, being
    uniformly packed in the sense of \cite{GoeTil:UPC}, as having minimum distance $5=2\times3-1$ and external distance $3$. In particular it is completely regular as per \cite[Thm 7]{GoeTil:UPC}. The fact that the error correcting capacity is $2$ implies that the first three degrees are $1,33=3\times 11,495={11 \choose 2}3^2$. The last degree follows by $2^{10}-1-33-495=495$.
    \qed\end{proof}
    More structural information on this graph is as follows.
    {\theorem The intersection array of the graph $\Gamma_{D^-}$ is $ \{33,30,15;1,2,15\}$.}
    \begin{proof}
    We sketch the three steps of the proof as follows. We know $\lambda$ and $\mu$ from Theorem \ref{lambmu}. The outer distribution matrix $B$ of an uniformly packed code is uniquely determined by $n$, $e$, $q$ and these two parameters \cite[Cor. 11.1.2]{Brouwer}. From the known intersection array of $H(n,4)$, and this data, the result follows by \cite[Th. 11.1.8]{Brouwer}, with $\Gamma=H(n,4)$, and $\Pi$ being the completely regular partition induced by the cosets of $D^-$.
    \qed\end{proof}

     {\theorem The automorphism group of the graph has order $2^{10}\cdot 54$, acts transitively on the vertices, and has two orbits on the edges,
     of size $2^9\cdot 6$ and $2^9\cdot 27$. The stabilizer of a vertex has a structure of type $(C_9 \rtimes C_3) \rtimes C_2$.
     The induced subgraph fixed by an element of the automorphism group, depending on the order $p$ of the element, can be the following, up to isomorphism:
     \begin{enumerate}
      \item[\rm (i)] the whole graph, \ $p=1;$
       \item[\rm (ii)] the null graph, \ $p=2;$
      \item[\rm (iii)] the disjoint union of $8$ complete graphs of order $4$, \ $p=2;$
      \item[\rm (iv)] the edge-free graph on $4$ vertices, \ $p=3;$
      \item[\rm (v)] the Hamming graph $H(2,4)$, \ $p=3$;
      \item[\rm (vi)] the edge-free graph on $2$ vertices, \ $p=6;$
      \item [\rm (vii)] the one-vertex graph, \ $p=9;$
      \item [\rm (viii)] the complete graph of order $4$, \ $p=9$.
     \end{enumerate}
     }
    \begin{proof}
     Since the graph is a Cayley graph, it is vertex transitive. The other assertions were established by Sage \cite{sage} computations.
    \qed\end{proof}

    The possible induced subgraphs fixed by an element of a prime order
    of the automorphism group of a distance-regular graph with intersection array $\{33,30,15;1,2,15\}$ were studied in~\cite{Makh:33-30-15-1-2-15}.
    The case (iv) of the theorem above is missing in the main Theorem of~\cite{Makh:33-30-15-1-2-15},
    which should be completed by the subcase where every two vertices of the
fixed subgraph are at distance three \cite{Makhnev:private2018}.

    The following strongly regular graph is found by standard spectral techniques, as indicated in \cite[p.\,428]{Brouwer}. For the next two results, we assume that the reader is familiar with the theory of duality in association schemes \cite{BanIt-1,Brouwer,Delsarte:1973}.
    {\theorem The graph $(\Gamma_{D^-})_2$ is strongly regular with parameters $(2^{10}, 495, 238, 240)$.}
\begin{proof}
The $P$-matrix of the association scheme underlying $\Gamma_{D^-}$ can be computed by the formulas in \cite[pp.\,135--136]{CaldGoeth:85:dual} or \cite[Th.\,5.25]{Delsarte:1973}, as
 $$\left(\begin{array}{cccc}
      1& 33& 495& 495 \\
  1 &  9& 15& -25\\
  1 &  1 &-17& 15 \\
 1 & -7&  15 & -9
  \end{array}\right). $$
  The second column is the spectrum of $(\Gamma_{D^-})_2$, discounting multiplicities. The result follows by the spectral characterization of SRGs \cite{BroHae:spectra}.
 \qed\end{proof}
    Strongly regular graphs with these parameters are known \cite{Brouwer:param}. We do not know if any of them is isomorphic to $(\Gamma_{D^-})_2$.

    There is a $Q$-analogue of the preceding results.
    {\theorem The Delsarte dual of the underlying association scheme of $\Gamma_{D^-}$ is $Q$-polynomial with Krein array $\{33,30,15;1,2,15\}$, multiplicities $1$, $33$, $495$, $495$ and valencies $1$, $198$, $495$, $330$. The second relation of that association scheme is a SRG of parameters $(2^{10}, 495,238, 240)$. }
    \begin{proof}
    Recall that the Delsarte dual of a coset graph is an association scheme on the codewords of the dual code \cite{Delsarte:1973}, called hereafter the {\em distance scheme}. If the weights are numbered $$w_1=6<w_2=8<w_3=10,$$ 
    then for $x,y\in {D^{-}}^\bot$, we define the relation $R_i$ as $x R_i y $ iff $\mathrm{wt}(x+y)=w_i$.
    By convention $x R_0 y$ iff $x=y$. This is exactly the situation of \cite{CaldGoeth:85:dual} with $e=2$, up to the fact that $D^-$ is not linear but only additive. By Delsarte duality, the multiplicities of this scheme are the degrees of
    $\Gamma_{D^-}$. Similarly, its valencies are the multiplicities of the spectrum of $\Gamma_{D^-}$. The Krein array coincides with the intersection array of the coset graph.
    The fact that the second relation of the distance scheme is a SRG comes from computation of the $P$-matrix which can be done as in \cite[pp. 135--136]{CaldGoeth:85:dual}.
    $$\left(\begin{array}{cccc}
      1& 198& 495& 330\\
  1 & 54 & 15& -70\\
  1 &  6 &-17&  10\\
 1 &-10&  15 & -6
  \end{array}\right).\eqno{\qed} $$
    \end{proof}
\begin{remark}
 By Delsarte duality, the P-matrix of the coset scheme is the Q-matrix of the distance scheme and conversely.
\end{remark}
\begin{remark}
 This Q-polynomial scheme appears as an open problem in the table of Jason Williford \cite{Williford:tbl}. We thank Bill Martin for pointing this out.
\end{remark}
\begin{remark}
 We do not know if the SRG constructed in that theorem is isomorphic with $(\Gamma_{D^-})_2$.
\end{remark}

     We also obtain a strongly regular graph on the cosets of $D^{3-}$.
    {\theorem  The graph $\Gamma_{D^{3-}}$ is a SRG of parameters $(64,27,10,12)$. The spectrum of $\Gamma_{D^{3-}}$ is $\{27^1,3^{36},(-5)^{27}\}$. }
\begin{proof}
The spectrum of $\Gamma_{D^{3-}}$ is easily computed from the weight distribution of Theorem \ref{distri2} upon applying \cite[11.1.11]{Brouwer}. The strong regularity follows then by the spectral characterization of SRG's \cite{BroHae:spectra}.
The parameters follow by the data of the spectrum \cite[Th. 9.1.3]{BroHae:spectra}.
\qed\end{proof}
        \begin{remark}
Depending on the punctured coordinates we obtain $3$ non-isomorphic SRG  $(64,27,10,12)$. One of them is isomorphic to the SRG corresponding to a linear $[9,3]$ two-weight quaternary code. All the codes $D^{3-}$ are, however, nonlinear.         
        \end{remark}
       \section{Conclusion and open problems}\label{s:concl}
       In this note, we have constructed the first additive non linear uniformly packed code in the history of the field. This solves a forty-four year old open problem of \cite{BZZ:1974:UPC}.
       The other values of $m$ in the Introduction are worth investigating, even if they lead to graph parameters beyond the tables of \cite{Brouwer,vDKT:DR} (the intersection array is of form $\{ 2^{2m+1}\!+1, 2^{2m+1}\!-2, 2^{2m}\!-1;1, 2, 2^{2m}\!-1\} $). The study of completely regular additive quaternary codes is only beginning.

       On another tack, the existence of an uniformly packed code in a Doob graph of diameter $11$ \cite[p.27]{Brouwer}, a distance-regular graph with the same parameters as the $H(11,4)$, is a goal worth pursuing. Note that perfect codes are known in Doob graphs \cite{Kro:perfect-doob}.

       That object would constitute a Galois ring analogue of the (punctured) dodecacode, and might lead to another distance-regular graph with intersection array $\{33,30,15;1,2,15\}$, non-isomorphic to the one considered in the current paper.

\begin{acknowledgements}
We thank Jack Koolen, Alexander Makhnev, and Bill Martin for helpful discussions and 
the anonymous referees for useful comments.
\end{acknowledgements}


\providecommand\href[2]{#2} \providecommand\url[1]{\href{#1}{#1}}
  \def\DOI#1{{\small {DOI}:
  \href{http://dx.doi.org/#1}{#1}}}\def\DOIURL#1#2{{\small{DOI}:
  \href{http://dx.doi.org/#2}{#1}}}

\end{document}